\documentclass[12pt,leqno]{article}
\usepackage{amssymb} 
\usepackage{amsmath} 
\usepackage{amsthm} 
\usepackage{graphicx} 
\usepackage{psfrag} 
\newtheorem{dummy}{}[section]

\newtheorem{theorem}[dummy]{Theorem}

\newtheorem{proposition}[dummy]{Proposition}
\newtheorem{lemma}[dummy]{Lemma}

\newtheorem{remark}[dummy]{Remark}

\begin{document}
\bibliographystyle{plain}
\title{ $\psi^{3}$ as an upper triangular matrix }
\author{Jonathan Barker and Victor Snaith}
\date{}
\maketitle
\begin{abstract}
In the $2$-local stable homotopy category the group of left-$bu$-module automorphisms of $bu \wedge bo$ which 
induce the identity on mod $2$ homology is isomorphic to the group of infinite upper triangular matrices with entries in the $2$-adic integers. 
We identify the conjugacy class of the matrix corresponding to $1 \wedge \psi^{3}$, where $\psi^{3}$ is the Adams operation.
\end{abstract}

\section{Introduction}

Let $bu$ and $bo$ denote the stable homotopy spectra representing $2$-adically completed 
unitary and orthogonal connective K-theory respectively. The main result of \cite{SnTri} is the existence of an isomorphism of groups
\[  \Psi :  U_{\infty}\mathbb{Z}_2   \stackrel{\cong}{\longrightarrow}  {\rm Aut}_{left-bu-mod}^{0}( bu \wedge bo )    .\]
Here $ U_{\infty}\mathbb{Z}_2 $ is the group of upper triangular matrices with coefficients in the $2$-adic integers and 
${\rm Aut}_{left-bu-mod}^{0}( bu \wedge bo )$ denotes the group of left $bu$-module automorphisms of $bu \wedge bo$ in the 
stable homotopy category of $2$-local spectra, which induce the identity on mod $2$ singular homology. The details of this isomorphism 
are recapitulated  in \S2.

This isomorphism is defined  up to inner automorphisms of $ U_{\infty}\mathbb{Z}_2 $. Given an important automorphism in 
${\rm Aut}_{left-bu-mod}^{0}( bu \wedge bo ) $ one is led to ask what is its conjugacy class in $ U_{\infty}\mathbb{Z}_2$.
By far the most important such automorphism is $1 \wedge \psi^{3}$, where $\psi^{3} : bo \longrightarrow bo$ denotes the Adams operation.

The following is our main result, which proved by combining the discussion of \S\ref{3.5} with Theorem \ref{4.2}.
\begin{theorem}{$_{}$}
\label{0.1}
 
Under the isomorphism $\Psi$ the automorphism $1 \wedge \psi^{3}$ corresponds to an element in the conjugacy class of the matrix
\[   \left( 
Ê\begin{array}{cccccc} 

1 & 1 & 0 & 0 & 0 & \ldots \\
\\
0 & 9 & 1 & 0 & 0 & \ldots \\
\\
0 & 0 & 9^{2} & 1 & 0 Ê& \ldots \\
\\
0 & 0 & 0 & 9^{3} & 1 Ê & \ldots \\
\\
\vdots & Ê\vdots & Ê\vdots & Ê\vdots & Ê\vdots & Ê\vdots 
\end{array} Ê\right).  \]
\end{theorem}

By techniques which are described in (\cite{A} pp.338-360) and reiterated in \S3, Theorem \ref{0.1} reduces to calculating 
the effect on $1 \wedge \psi^{3}$ on $\pi_{*}(bu \wedge bo)$ modulo torsion. The difficulty arises  because, in order to identify
$\Psi^{-1}( 1 \wedge \psi^{3})$ one must compute the map on homotopy modulo torsion in terms of an unknown $2$-adic basis defined 
in terms of the Mahowald splitting of $bu \wedge bo$ (see \S2 and \S3). On the other hand a very convenient $2$-adic basis is defined
in \cite{ccw} and the crucial fact is that $1 \wedge \psi^{3}$ acts on the second basis by the matrix of Theorem \ref{0.1}. This fact was 
pointed out to one of us (VPS) by Francis Clarke in 2001 and led to the confident prediction appearing as a footnote in (\cite{SnTri} p.1273).
Verifying the prediction has proved a little more difficult than first imagined!

Once one has Theorem \ref{0.1} a number of homotopy problems become merely a matter of matrix algebra. In \S5 we give an example 
concerning the maps $1 \wedge (\psi^{3} - 1) (\psi^{3} - 9) \ldots (\psi^{3} - 9^{n-1})$ where we prove a vanishing result
(Theorem \ref{5.2}) which is closely related to the main theorem of \cite{Mg}, as explained in Remark \ref{5.3}. In subsequent papers we shall give further
applications to homotopy theory and to algebraic K-theory.

\section{$2$-adic homotopy of $bu \wedge bo$}
\begin{dummy}
\label{2.1}
\begin{em}

Let $bu$ and $bo$ denote the stable homotopy spectra representing $2$-adically completed 
unitary and orthogonal connective K-theory respectively.
We shall begin by recalling the $2$-local homotopy decomposition of $bu \wedge bo$ which is one of a number of similar 
results which were discovered by Mark Mahowald in the 1970's. These results may be proved in several ways \cite{A}, \cite{Mah} and \cite{Mg}.
For notational reasons we shall refer to the proof, a mild modification of (\cite{A} pp.190-196), which appears in (\cite{SnTri} \S2).

Consider the second loopspace of the $3$-sphere, $\Omega^{2} S^{3}$.
There exists a model for $\Omega^{2} S^{3}$ which is filtered by
finite complexes (\cite{BP} , \cite{Snspl})
\[ S^{1} = F_{1} \subset F_{2} \subset F_{3} \subset  \ldots  \subset 
 \Omega^{2} S^{3} = \bigcup_{k \geq 1}  F_{k} \]
and there is a stable homotopy equivalence, an example of the so-called Snaith splitting, of the form
\[ \Omega^{2} S^{3}  \simeq   \vee_{k \geq 1}  F_{k}/F_{k-1} .  \]
There is a $2$-local homotopy equivalence of left-$bu$-module spectra (see \cite{SnTri} Theorem 2.3(ii)) of the form
\[   \hat{L} :  \vee_{k \geq 0}  bu \wedge  (  F_{4k}/F_{4k-1})  
\stackrel{\simeq}{\longrightarrow }  bu  \wedge  bo.  \]
The important fact about this homotopy equivalence is that its induced map on mod $2$ 
homology is a specific isomorphism which is described in (\cite{SnTri} \S2.2) 

From this decomposition we obtain left-$bu$-module spectrum maps of the form
\[   \iota_{k,l} :  bu \wedge (F_{4k}/F_{4k-1})  \longrightarrow  bu \wedge (F_{4l}/F_{4l-1}) \]
where $\iota_{k,k} = 1$, $\iota_{k,l} = 0$ if $l > k$ and, as explained in (\cite{SnTri} \S3.1),  
$\iota_{k,l}$ is defined up to multiplication by a $2$-adic unit when $k>l$.

Consider the ring of left
$bu$-module endomorphisms of degree zero in the stable homotopy category of spectra \cite{A}, 
which we shall denote by 
\linebreak
${\rm End}_{left-bu-mod}( bu \wedge bo )$. The group of units in this ring
will be denoted by ${\rm Aut}_{left-bu-mod}( bu \wedge bo )$, the group of homotopy classes of 
left $bu$-module homotopy equvialences and let ${\rm Aut}_{left-bu-mod}^{0}( bu \wedge bo )$
denote the subgroup of left $bu$-module homotopy equivalences which induce the identity map on 
$H_{*}(bu \wedge bo ; \mathbb{Z}/2)$.

Let $U_{\infty}\mathbb{Z}_2$ denote the group of infinite, invertible upper triangular matrices with entries
in the $2$-adic integers. That is, $X = (X_{i,j}) \in U_{\infty}\mathbb{Z}_2$ if $X_{i,j} \in \mathbb{Z}_2$
for each pair of integers $0 \leq i,j$ and $X_{i,j} = 0$ if $j > i$ and $X_{i,i}$ is a $2$-adic unit. This upper
trianglular group is {\em not} equal to the direct limit $ \lim_{\stackrel{\rightarrow}{n}}  \ U_{n}\mathbb{Z}_2 $
of the finite upper triangular groups.
The main result of \cite{SnTri} is the existence of an isomorphism of groups
\[  \Psi :  U_{\infty}\mathbb{Z}_2   \stackrel{\cong}{\longrightarrow}  {\rm Aut}_{left-bu-mod}^{0}( bu \wedge bo )    .\]
By the Mahowald decomposition of $bu \wedge bo$ the existence of $\psi$ is equivalent to an isomorphism of the form
\[ \Psi :  U_{\infty}\mathbb{Z}_2  \stackrel{\cong}{\longrightarrow}  
{\rm Aut}_{left-bu-mod}^{0}( \vee_{k \geq 0}  bu \wedge  ( F_{4k}/F_{4k-1})  ) . \]
If we choose $\iota_{k,l} $ to satisfy 
$ \iota_{k,l} =  \iota_{l+1,l} \iota_{l+2,l+1} \ldots \iota_{k,k-1}  $
for all $k - l \geq 2$ then, for $X \in  U_{\infty}\mathbb{Z}_{2}$, we define (\cite{SnTri} \S3.2)
\[ \Psi(X^{-1}) = \sum_{l \leq k}    X_{l,k} \iota_{k,l} : 
 bu \wedge ( \vee_{k \geq 0}  F_{4k}/F_{4k-1} )  \longrightarrow 
 bu \wedge ( \vee_{k \geq 0}  F_{4k}/F_{4k-1} ).   \]
The ambiguity in the definition of the $\iota_{k,l}$'s implies that $\Psi$ is defined up to conjugation 
by a diagonal matrix in $U_{\infty}\mathbb{Z}_2$.

\end{em}
\end{dummy}
\begin{dummy}{Bases for $\frac{\pi_*(bu\wedge
bo)\otimes\mathbb{Z}_2}{{\rm Torsion}}$}
\label{2.2}
\begin{em}

Let $G_{s,t}$ denote the $2$-adic homotopy group modulo torsion
$$G_{s,t} = \frac{\pi_{s}(bu\wedge F_{4t}/F_{4t-1})\otimes\mathbb{Z}_2}{{\rm Torsion}}$$
so 
\[  G_{*,*} =  \oplus_{s,t} \ \frac{\pi_{s}(bu\wedge F_{4t}/F_{4t-1})\otimes\mathbb{Z}_2}{{\rm Torsion}}
\cong   \frac{\pi_{*}(bu\wedge
bo)\otimes\mathbb{Z}_2}{{\rm Torsion}}.   \]
From \cite{A} or \cite{SnTri} 
\[  G_{s,t} \cong \left\{  
\begin{array}{ll}
\mathbb{Z}_{2}  & {\rm if}  \  s \  {\rm even}, \ s \geq 4t,  \\
\\
0 & {\rm otherwise} 
\end{array} \right.   \]
and if $\tilde{G}_{s,t} $ denotes $\pi_{s}(bu\wedge F_{4t}/F_{4t-1})\otimes\mathbb{Z}_2$ then
$\tilde{G}_{s,t} \cong G_{s,t} \oplus W_{s,t}$ where $W_{s,t}$ is a finite, elementary abelian $2$-group.

In \cite{ccw} a $\mathbb{Z}_{2}$-basis is given for $G_{*,*}$ consisting of elements lying in the subring
$\mathbb{Z}_2[u/2,v^2/4]$ of $\mathbb{Q}_2[u/2,v^2/4]$. One starts with the elements
\[  c_{4k}=\Pi_{i=1}^k\Bigg(\frac{v^2-9^{i-1}u^2}{9^k-9^{i-1}}\Bigg),\hspace{20pt}k=1,2,\ldots  \]
and ``rationalises'' them, after the manner of (\cite{A} p.358), to obtain elements of
$\mathbb{Z}_2[u/2,v^2/4]$. In order to describe this basis we shall require a few well-known preparatory results
about $2$-adic valuations.
\end{em}
\end{dummy}
\begin{proposition}{$_{}$}
\label{2.3}

 For any integer $n \geq 0$, $9^{2^n}-1=2^{n+3}(2s+1)$ for some $s\in\mathbb{Z}$.
\end{proposition}
\vspace{2pt}

{\bf Proof}
\vspace{2pt}

We prove this by induction on $n$, starting with $9-1=2^3 $. Assuming the result is true for $n$, we have
\[ \begin{array}{ll}
 9^{2^{(n+1)}}-1 & = (9^{2^n}-1)(9^{2^n}+1) \\ 
\\
&=  (9^{2^n}-1)(9^{2^n}-1+2)   \\
\\
&=  2^{n+3}(2s+1)(2^{n+3}(2s+1)+2) \\
\\ 
&= 2^{n+4}(2s+1)\underbrace{(2^{n+2}(2s+1)+1)}_{odd}
\end{array}  \] 
as required. $\Box$

\begin{proposition}{$_{}$}
\label{2.4} 

For any integer $l \geq 0$, $9^l-1=2^{\nu_2(l)+3}(2s+1)$ for some $s\in\mathbb{Z}$, where $\nu_2(l)$ denotes the
$2$-adic valuation of $l$.
\end{proposition}
\vspace{2pt}

{\bf Proof}
\vspace{2pt} 

Write $l = 2^{e_{1}} + 2^{e_{2}} + \ldots  + 2^{e_{k}}$ with $0 \leq e_{1} < e_{2} < \ldots < e_{k}$ so that
$\nu_2(l) = e_{1}$. Then, by Proposition \ref{2.3},
\[ \begin{array}{ll}
9^l-1  &  =  9^{2^{e_1}+2^{e_2}+\cdots+2^{2_k}}-1 \\
\\
& = ((2s_{1} + 1) 2^{e_{1} + 3} + 1))  \ldots ((2s_{k} + 1) 2^{e_{k} + 3} + 1))  - 1  \\
\\
& \equiv   (2s_{1} + 1) 2^{e_{1} + 3}  \  ({\rm modulo}  2^{e_{1} + 4 }) \\
\\
& = 2^{e_{1} + 3} (2t+1)  
\end{array} \]
as required. $\Box$

\begin{proposition}{$_{}$}

\label{2.5} 

For any integer $l \geq 1$, $\prod_{i=1}^{l} \ (9^{l} - 9^{i-1}) = 2^{\nu_{2}(l!) + 3l}(2s+1)$ for some $s\in\mathbb{Z}$.
\end{proposition}
\vspace{2pt}

{\bf Proof}
\vspace{2pt} 

By Proposition \ref{2.4} we have
\[ \begin{array}{ll}
\prod_{i=1}^{l} \ (9^{l} - 9^{i-1})  & =  \prod_{i=1}^{l} \ (9^{l-i+1} - 1) 9^{i-1}  \\
\\
& =  \prod_{i=1}^{l} \ 2^{\nu_{2}(l-i+1) + 3} (2t_{i}+1)  9^{i-1} \\
\\
& = (2t+1)  2^{\nu_{2}(l!) + 3l},
\end{array} \]
as required. $\Box$
\begin{proposition}
\label{2.6} 

For any integer $l \geq 0$, $ 2^{\nu_{2}(l!) + 3l} =  2^{4l - \alpha(l)}$ where $\alpha(l)$ 
is equal to the number of $1$'s in the dyadic expansion of $l$.
\end{proposition}
\vspace{2pt}

{\bf Proof}
\vspace{2pt} 

Write $l = 2^{e_{1}} + 2^{e_{2}} + \ldots  + 2^{e_{k}}$ with $0 \leq e_{1} < e_{2} < \ldots < e_{k}$
so that $\alpha(l) = k$.

$9^l-1=2^{\nu_2(l)+3}(2s+1)$ for some $s\in\mathbb{Z}$, where $\nu_2(l)$ denotes the
$2$-adic valuation of $l$.
Then
\[  \begin{array}{crclllllll}
\nu_2(l!)&=&&2^{\alpha_1-1}&+&2^{\alpha_2-1}&+&\cdots&+&2^{\alpha_k-1}\\
&&+&2^{\alpha_1-2}&+&2^{\alpha_2-2}&+&\cdots&+&2^{\alpha_k-2}\\ &&&\vdots&&&&&&\vdots\\
&&+&1&+&2^{\alpha_2-\alpha_1}&+&\cdots&+&2^{\alpha_k-\alpha_1}\\ &&&&+&1&+&\cdots&+&2^{\alpha_k-\alpha_2}\\ &&&&&&&&+&1 
\end{array}  \]
because the first row counts the multiples of $2$ less than or equal to $l$, the second row counts the multiples of $4$,
the third row counts multiples of $8$ and so on.
Adding by columns we obtain
\[  \nu_2(l!) = 2^{\alpha_1}-1+2^{\alpha_2}-1+\cdots+2^{\alpha_k}-1\\ =  l - k \]
which implies that $2^{3l+\nu_2(l!)}=2^{3l+l-\alpha(l)}=2^{4l-\alpha(l)}$, as required. $\Box$
\begin{dummy}{Bases continued}
\label{2.7}
\begin{em}

Consider the elements 
$c_{4k}=\Pi_{i=1}^k  \big(\frac{v^2-9^{i-1}u^2}{9^k-9^{i-1}}\big)$, introduced in \S\ref{2.2}, 
for a particular $k = 1, 2, \ldots$. 
For completeness write $c_{0} = 1$ so that $c_{4k} \in \mathbb{Q}_2[u/2,v^2/4]$.
Since the degree of the numerator of $c_{4k}$ is $2k$, Proposition \ref{2.6} implies that
\[ f_{4k}=2^{4k-\alpha(k)-2k}c_{4k}=2^{2k-\alpha(k)}\Pi_{i=1}^k\big(\frac{v^2-9^{i-1}u^2}{9^k-9^{i-1}}\big) \]
lies in  $\mathbb{Z}_2[u/2,v^2/4]$ but $2^{4k-\alpha(k)-2k-1}c_{4k} \notin \mathbb{Z}_2[u/2,v^2/4]$.
Similarly $(u/2) f_{4k} = 2^{4k-\alpha(k)-2k-1} u c_{4k} \in \mathbb{Z}_2[u/2,v^2/4]$ but
$2^{4k-\alpha(k)-2k-2} u c_{4k} \notin \mathbb{Z}_2[u/2,v^2/4]$ and so on. This process is the 
``rationalisation yoga'' referred to in \S\ref{2.2}. One forms $u^{j} c_{4k}$ and then multiplies by the smallest 
positive power of $2$ to obtain an element of $\mathbb{Z}_2[u/2,v^2/4]$. 

By Proposition \ref{2.6}, starting with $f_{4l}  =2^{4l-\alpha(l)-2l}c_{4l}$ this process produces the
following set of
elements of $ \mathbb{Z}_2[u/2,v^2/4]$
\[ \begin{array}{l}
f_{4l}, \ (u/2)f_{4l}, \ (u/2)^{2}f_{4l},  \  \ldots  \  , \  (u/2)^{2l - \alpha(l)} f_{4l}, \\
\\
\hspace{50pt}  u (u/2)^{2l - \alpha(l)} f_{4l},  \
 u^{2} (u/2)^{2l - \alpha(l)} f_{4l},  \  u^{3} (u/2)^{2l - \alpha(l)} f_{4l},  \   \ldots .
\end{array}  \]

As explained in (\cite{A} p.352 et seq), the Hurewicz homorphism defines an injection of graded groups of the form
\[  \frac{\pi_{*}(bu \wedge bo) \otimes \mathbb{Z}_{2}}{{\rm Torsion}} \longrightarrow  \mathbb{Q}_2[u/2,v^2/4] \]
which, by the main theorem of \cite{ccw}, induces an isomorphism between $\frac{\pi_{*}(bu \wedge bo) \otimes \mathbb{Z}_{2}}{{\rm
Torsion}}$ and the free graded $\mathbb{Z}_{2}$-module whose basis consists of the elements of 
\linebreak
$ \mathbb{Z}_2[u/2,v^2/4]$
listed above for $l = 0, 1, 2, 3, \ldots$.

From this list we shall be particularly interested in the elements whose degree is a multiple of $4$. Therefore
denote by $g_{4m,4l} \in \mathbb{Z}_2[u/2,v^2/4]$ for $l
\leq m$ the element produced from $f_{4l}$ in degree $4m$. Hence, for $m \geq l$, $g_{4m,4l}$ is given by the formula
\[ g_{4m,4l} =    \left\{
\begin{array}{ll}
u^{2m-4l+\alpha(l)}\big[\frac{u^{2l-\alpha(l)}f_{4l}}{2^{2l-\alpha(l)}}\big] &   {\rm if}  \   4l-\alpha(l) \leq 2m ,  \\
\\
\big[\frac{u^{2(m-l)f_{4l}}}{2^{2(m-l)}}\big]  &  {\rm if}  \   4l-\alpha(l) > 2m   .

\end{array} \right.  \]
\end{em}
\end{dummy}
\begin{lemma}{$_{}$}
\label{2.8}

In the notation of \S\ref{2.2}, let $\Pi$ denote the projection
\[ \Pi :  \frac{\pi_{*}(bu\wedge
bo)\otimes\mathbb{Z}_2}{{\rm Torsion}}  \cong  G_{*,*}  \longrightarrow   G_{*,k} = \oplus_{m} \  G_{m,k}  . \]
Then $\Pi(g_{4k,4i}) = 0$ for all $i < k$.
\end{lemma}
\vspace{2pt}

{\bf Proof}
\vspace{2pt}

Since $G_{m,k}$ is torsion free it suffices to show that  $\Pi(g_{4k,4i})$ vanishes in $ G_{*,k}  \otimes \mathbb{Q}_{2}$.
When $i < k$, by definition
\[  g_{4k,4i} \in   u^{2k - 2i} \frac{\pi_{4i}(bu\wedge
bo)\otimes\mathbb{Z}_2}{{\rm Torsion}} \otimes \mathbb{Q}_{2}   \subset   \frac{\pi_{4k}(bu\wedge
bo)\otimes\mathbb{Z}_2}{{\rm Torsion}} \otimes \mathbb{Q}_{2} . \]
However $\Pi$ projects onto $ \oplus_{s} \ \frac{\pi_{s}(bu\wedge F_{4k}/F_{4k-1})\otimes\mathbb{Z}_2}{{\rm Torsion}}$
and commutes with multiplication by $u$ so the result follows from the fact that the homotopy of 
$bu\wedge F_{4k}/F_{4k-1}$ is trivial in degrees less than $4k$ (see \cite{SnTri} \S3). $\Box$
\begin{dummy}
\label{2.9}
\begin{em}

Recall from \S\ref{2.2} that $G_{4k,k} \cong \mathbb{Z}_{2}$ for $k = 0,1,2,3, \ldots $ so we may choose a generator $z_{4k}$ 
for this group as a module over the $2$-adic integers (with the convention that $z_{0} = f_{0} = 1$). Let
$\tilde{z}_{4k}$ be any choice of an element in the $2$-adic homotopy group $\tilde{G}_{4k,k} \cong G_{4k,k} \oplus W_{4k,k}$
whose first coordinate is $z_{4k}$.
\end{em}
\end{dummy}
\begin{lemma}{$_{}$}
\label{2.10}

Let $B$ denote the exterior subalgebra of $\mathbb{Z}/2$ Steenrod algebra generated by $Sq^{1}$ and $Sq^{0,1}$. 
In the collapsed Adams spectral sequence (see \cite{A} or \cite{SnTri})
\[  \begin{array}{l}
E_{2}^{s,t} \cong Ext_{B}^{s,t}(  H^{*}(F_{4k}/F_{4k-1} ; \mathbb{Z}/2) , \mathbb{Z}/2 )   \\
\\
\hspace{50pt}   \Longrightarrow
\pi_{t-s}(  bu \wedge (F_{4k}/F_{4k-1})  ) \otimes  \mathbb{Z}_{2}  
\end{array}  \]
the homotopy class $\tilde{z}_{4k}$ is represented either in $E_{2}^{0,4k}$ or $E_{2}^{1.4k+1}$.
\end{lemma}
\vspace{2pt}

{\bf Proof}
\vspace{2pt}

Recall from \S\ref{2.2} that $\pi_{4k}(  bu \wedge (F_{4k}/F_{4k-1})  ) \otimes  \mathbb{Z}_{2} = \tilde{G}_{4k,k} \cong
\mathbb{Z}_{2} \oplus W_{4k,k}$. The following behaviour of the filtration coming from the spectral sequence is well-known, being
explained in \cite{A}. The group $\tilde{G}_{4k,k}$ has a filtration
\[ \ldots  \subset  F^{i} \subset \ldots  F^{2} \subset F^{1} \subseteq F^{0} = \tilde{G}_{4k,k}  \]
with $F^{i}/F^{i+1} \cong E_{2}^{i, 4k+i}$ and $2 F^{i} \subseteq F^{i+1}$. Also $2 \cdot   W_{4k,k} = 0$, every non-trivial
element of $ W_{4k,k}$ being represented in $E_{2}^{0,4k}$ . Furthermore for $i= 1, 2, 3, \ldots$ we have
$2 F^{i} = F^{i+1}$ and $F^{1} \cong \mathbb{Z}_{2}$.

Now suppose that $\tilde{z}_{4k}$ is represented in $E_{2}^{j,4k+j}$ for $j \geq 2$ then $\tilde{z}_{4k} \in F^{j}$. From the
multiplicative structure of the spectral sequence there exists a generator $\hat{z}_{4k}$ of $F^{1}$ such that
$2^{j} \hat{z}_{4k}$ generates $F^{j+1}$ and therefore $2^{j} \gamma \hat{z}_{4k} = 2 \tilde{z}_{4k}$ for
some $2$-adic integer $\gamma$. Hence $2( 2^{j-1} \gamma \hat{z}_{4k} - \tilde{z}_{4k}) = 0$ and so
$ 2^{j-1} \gamma \hat{z}_{4k} - \tilde{z}_{4k}  \in W_{4k,k}$ which implies  the contradiction that the generator $z_{4k}$ is
divisible by
$2$ in
$G_{4k,k}$. $\Box$ 
\begin{theorem}{$_{}$}
\label{2.11}

In the notation of \S\ref{2.7} and \S\ref{2.9}
\[  z_{4k} = \Sigma_{i=0}^k2^{\beta(k,i)}  \lambda_{4k,4i} g_{4k,4i} \in  
\frac{\pi_{4k}(bu \wedge bo) \otimes \mathbb{Z}_{2}}{{\rm Torsion}}  \] 
with $\lambda_{s,t}\in\mathbb{Z}_2$, $\lambda_{4k,4k}\in\mathbb{Z}_2^*$ and
\[  \beta(k,i) =  \left\{  
\begin{array}{ll}
  4(k- i) - \alpha(k) + \alpha(i)   &   {\rm if}  \  4i - \alpha(i) > 2k ,  \\
\\
  2k - \alpha(k) & {\rm if}  \   4i - \alpha(i) \leq 2k.
\end{array} \right.  \]
\end{theorem} 
\vspace{2pt}

{\bf Proof}
\vspace{2pt}

From \cite{ccw}, as explained in \S\ref{2.7}, a $\mathbb{Z}_{2}$-module basis for $G_{4k,*}$ is given by
$  \{  g_{4k,4l}  \}_{ 0  \leq l  \leq   k}$.
Hence there is a relation of the form
\[  z_{4k} = \lambda_{4k,4k} g_{4k,4k} +  \tilde{\lambda}_{4k,4(k-1)} g_{4k,4(k-1)}  +  \ldots  +
\tilde{\lambda}_{4k,0}g_{4k,0}  \]
where $\tilde{\lambda}_{4k,4i}$ and $\lambda_{4k,4k}$ are $2$-adic integers.   
Applying the projection 
\linebreak
$\Pi:G_{4k,*} \longrightarrow G_{4k,k}$ we see that
$z_{4k} = \Pi(z_{4k}) =  \lambda_{4k,4k} \Pi( g_{4k,4k})$, by Lemma \ref{2.8}. Hence, if $ \lambda_{4k,4k}$ is not a $2$-adic unit,
then $z_{4k}$ would be divisible by $2$ in $G_{4k,k}$ and this is impossible since $z_{4k}$ is a generator, by definition.
 
Multiplying the relation
\[ z_{4k} =  \lambda_{4k,4k} \Pi( g_{4k,4k})  =   \lambda_{4k,4k} \Pi( f_{4k})  \in G_{4k,k} . \]
 by $(u/2)^{ 2k - \alpha(k)}$ we obtain $(u/2)^{ 2k - \alpha(k)} z_{4k} =  \lambda_{4k,4k} \Pi( (u/2)^{ 2k - \alpha(k)}f_{4k})$,
which lies in $G_{8k - 2 \alpha(k),k}$, by the discussion of \S\ref{2.7}. Therefore, in $G_{8k - 2 \alpha(k),k} \otimes
\mathbb{Q}_{2}$ we have the relation
\[     (u / 2)^{2k-\alpha(k)}  z_{4k}
=   (u / 2)^{2k-\alpha(k)}  f_{4k}   +  \sum_{i=0}^{k-1}    \tilde{\lambda}_{4k, 4i}   (u / 2)^{2k-\alpha(k)} g_{4k,4i}. \]
Since the left hand side of the equation lies in $G_{8k - 2 \alpha(k),k}$, the $\mathbb{Q}_{2}$ coefficients must all be $2$-adic
integers once we re-write the right hand side in terms of the basis of \S\ref{2.7}.

For $i = 0, 1, \ldots , k-1$ 
\[  \begin{array}{ll}
(u / 2)^{2k-\alpha(k)} g_{4k,4i} & = \left\{
\begin{array}{ll}
  \frac{u^{2k-\alpha(k) + 2k - 4i + \alpha(i) + 2i-\alpha(i)}}{2^{2k-\alpha(k) + 2i-\alpha(i)}} f_{4i}  &  {\rm if}  \  4i -
\alpha(i) \leq 2k, 
\\
\\
 \frac{u^{2k-\alpha(k) + 2k- 2i}}{2^{2k-\alpha(k) + 2k- 2i}}  f_{4i}   &  {\rm if}  \  4i - \alpha(i)  >   2k
\end{array} \right.   \\
\\
&   = \left\{
\begin{array}{ll}
  \frac{u^{4k - 2i -\alpha(k)  }}{2^{2k + 2i - \alpha(k) - \alpha(i)}} f_{4i}  &  {\rm if}  \  4i -
\alpha(i) \leq 2k, 
\\
\\
 \frac{u^{4k - 2i -  \alpha(k) }}{2^{4k -  2i -\alpha(k) }}  f_{4i}   &  {\rm if}  \  4i - \alpha(i)  >   2k.
\end{array} \right.  
\end{array}  \]
Now we shall write $(u / 2)^{2k-\alpha(k)} g_{4k,4i}$ as a power of $2$ times a generator derived from $f_{4i}$ 
in \S\ref{2.7} (since we did not define any generators called $g_{4k+2,4i}$ the generator in question will be
$g_{8k - 2 \alpha(k),4i}$ only when $\alpha(k)$ is even).

Assume that $ 4i - \alpha(i) \leq 2k$ so that $2i - \alpha(i) \leq 4k - 2i - \alpha(k)$ and 
\[   \frac{u^{4k - 2i -\alpha(k)  }}{2^{2k + 2i - \alpha(k) - \alpha(i)}} f_{4i} =  
 \frac{1}{2^{2k  - \alpha(k) }}u^{4k - 4i -\alpha(k) + \alpha(i) } (u/2)^{2i - \alpha(i)}
f_{4i} \]
which implies that $   \tilde{\lambda}_{4k, 4i} $ is divisible by $2^{2k  - \alpha(k) }$ in the $2$-adic integers, as required.

Finally assume that $4i - \alpha(i) > 2k$. We have  
$2i - \alpha(i) \leq 4k - 2i - \alpha(k)$ also.
To see this observe that $\alpha(i) + \alpha(k-i) - \alpha(k) \geq 0$ because, by Proposition \ref{2.6}, this equals the $2$-adic
valuation of the binomial coefficient $( \stackrel{k}{i} )$. Therefore
\[  \alpha(k) - \alpha(i) \leq \alpha(k-i) \leq  k-i < 4(k-i). \]
Then, as before,
\[   \frac{u^{4k - 2i -  \alpha(k) }}{2^{4k -  2i -\alpha(k) }}  f_{4i} =  
 \frac{1}{2^{4k - 4i - \alpha(k) + \alpha(i)}}u^{4k - 4i -\alpha(k) + \alpha(i) } (u/2)^{2i - \alpha(i)}
f_{4i} \]
which implies that $   \tilde{\lambda}_{4k, 4i} $ is divisible by $2^{4k - 4i - \alpha(k) + \alpha(i)}$ in the 
$2$-adic integers, as required.
$\Box$
\begin{theorem}{$_{}$}
\label{2.12}

(i)  \  In the collapsed Adams spectral sequence and the notation of Lemma \ref{2.10} $\tilde{z}_{4k}$ may be 
chosen to be represented in $E_{2}^{0,4k}$.

(ii)  \  In fact, $\tilde{z}_{4k}$ may be taken to be the smash product of the unit $\eta$ of the $bu$-spectrum
with the inclusion of the bottom cell $j_{k}$ into $ F_{4k}/F_{4k-1}$
\[ S^{0} \wedge S^{4k}   \stackrel{ \eta  \wedge  j_{k}}{\longrightarrow}  bu \wedge  F_{4k}/F_{4k-1} . \]
\end{theorem}
\vspace{2pt}

{\bf Proof}
\vspace{2pt}

For part (i), suppose that $\tilde{z}_{4k} $ is represented in $E_{2}^{1,4k+1}$. By Lemma \ref{2.10} we must show that this leads to
a contradiction. From \cite{SnTri} we know that on the $s=1$ line the non-trivial groups are precisely
$E_{2}^{1,4k+1}, E_{2}^{1,4k+3}, \ldots , E_{2}^{1,8k + 2 -  2 \alpha(k)}$ which are all of order two.
From the multiplicative structure of the spectral sequence, if a homotopy class $w$ is represented $E_{2}^{j, 4k+2j-1}$ and 
$E_{2}^{j, 4k+2j+1}$ is non-zero then there is a homotopy class $w'$ represented in $E_{2}^{j, 4k+2j+1}$ 
such that $2w' = uw$. Applied to $\tilde{z}_{4k}$ this implies that the homotopy element
$ u^{2k - \alpha(k) + 1} \tilde{z}_{4k}$ is divisible by $ 2^{2k - \alpha(k) + 1}$. Hence 
$ u^{2k - \alpha(k) + 1} z_{4k}$ is divisible by $ 2^{2k - \alpha(k) + 1}$ in $G_{*,*}$, which contradicts the proof of Theorem
\ref{2.11}. 

For part (ii) consider the Adams spectral sequence
\[  \begin{array}{l}
E_{2}^{s,t} = Ext_{B}^{s,t}(  H^{*}(F_{4k}/F_{4k-1} ; \mathbb{Z}/2)  , \mathbb{Z}/2 )   
  \Longrightarrow
\pi_{t-s}(  bu  \wedge F_{4k}/F_{4k-1} ) \otimes  \mathbb{Z}_{2} . 
\end{array}  \]
We have an isomorphism 
\[ E_{2}^{0,t} = Hom(  \frac{  H^{t}(F_{4k}/F_{4k-1} ; \mathbb{Z}/2) }{Sq^{1} H^{t-1}(F_{4k}/F_{4k-1} ; \mathbb{Z}/2)
+ Sq^{0,1} H^{t-3}(F_{4k}/F_{4k-1} ; \mathbb{Z}/2)}  ,  \mathbb{Z}/2) . \]
The discussion of the homology groups $ H_{*}(F_{4k}/F_{4k-1} ; \mathbb{Z}/2) $ given in (\cite{A} p.341; see also \S\ref{3.1})
shows that $E_{2}^{0,4k} \cong \mathbb{Z}/2$ generated by the Hurewicz image of $\eta  \wedge  j_{k}$. Therefore the generator of
$E_{2}^{0,4k}$ represents $\eta  \wedge  j_{k}$. Since there is only one non-zero element in $E_{2}^{0,4k}$ it must also
represent $\tilde{z}_{4k}$, by part (i), which completes the proof.
$\Box$

\section{The Matrix}
\begin{dummy}
\label{3.1}
\begin{em}

Consider the left-$bu$-module spectrum map of \S\ref{2.1}
\[   \iota_{k,l} :  bu \wedge (F_{4k}/F_{4k-1})  \longrightarrow  bu \wedge (F_{4l}/F_{4l-1}) \]
when $l > k$. This map is determined up to homotopy by its restriction, via the unit of $bu$, to $(F_{4k}/F_{4k-1})$.
By S-duality this restriction is equivalent to a map of the form
\[ S^{0} \longrightarrow  D(F_{4k}/F_{4k-1})  \wedge  bu \wedge (F_{4l}/F_{4l-1}) ,  \]
which $DX$ denotes the S-dual of $X$. Maps of this form are studied by means of the (collapsed) Adams 
spectral sequence (see \cite{SnTri} \S3.1), where $B$ is as in Lemma \ref{2.10},
\[  \begin{array}{l}
E_{2}^{s,t} = Ext_{B}^{s,t}(  H^{*}(D(F_{4k}/F_{4k-1}) ; \mathbb{Z}/2) 
\otimes H^{*}(F_{4l}/F_{4l-1} ; \mathbb{Z}/2) , \mathbb{Z}/2 )   \\
\\
\hspace{50pt}   \Longrightarrow
\pi_{t-s}(D(F_{4k}/F_{4k-1})  \wedge
(F_{4l}/F_{4l-1})  \wedge bu) \otimes  \mathbb{Z}_{2} . 
\end{array}  \]

Recall from (\cite{A} p.332) that $\Sigma^{a}$ is the (invertible) $B$-module given by $\mathbb{Z}/2$
in degree $a$, $\Sigma^{-a} = Hom( \Sigma^{a} , \mathbb{Z}/2)$ and $I$ is the augmentation ideal,
$I = ker( \epsilon : B \longrightarrow  \mathbb{Z}/2 )$. Hence, if $b >0$,
$I^{-b} = Hom( I^{b} , \mathbb{Z}/2)$, where $I^{b}$ is the $b$-fold tensor product of $I$. 
These duality identifications may be verified using the criteria of (\cite{A} p.334 Theorem 16.3)
for identifying $ \Sigma^{a} I^{b}$.

In (\cite{A} p.341) it is shown that the $B$-module given by
\[  H^{-*}(D(F_{4k}/F_{4k-1}) ; \mathbb{Z}/2) 
\cong   H_{*}(F_{4k}/F_{4k-1} ; \mathbb{Z}/2)   \]
is stably equivalent to $\Sigma^{2^{r-1}+1} I^{2^{r-1}-1}$ when $0 < 4k = 2^{r}$. Therefore
\linebreak
$ H^{*}(D(F_{4k}/F_{4k-1}) ; \mathbb{Z}/2)$ is stably equivalent to 
$\Sigma^{-(2^{r-1}+1)} I^{1 - 2^{r-1}}$ when $0 < 4k = 2^{r}$.
If $k$ is not a power of two we may write $4k = 2^{r_{1}} + 2^{r_{2}} + \ldots + 2^{r_{t}}$ with
$2 \leq r_{1} < r_{2} < \ldots < r_{t}$. In this case
\[  H_{*}(F_{4k}/F_{4k-1} ; \mathbb{Z}/2)  \cong  \otimes_{j = r_{1}}^{r_{t}}  \  
 H_{*}(F_{2^{j}}/F_{2^{j}-1} ; \mathbb{Z}/2)   \]
which is stably equivalent to $ \Sigma^{ 2k + \alpha(k)}
I^{ 2k - \alpha(k) }$,
where $\alpha(k)$ equals the number of $1$'s in the dyadic expansion of $k$, as in Proposition \ref{2.6}. Similarly,
$ H^{*}(D(F_{4k}/F_{4k-1}) ; \mathbb{Z}/2)$ is stably equivalent to $ \Sigma^{ -2k - \alpha(k)}
I^{   \alpha(k) - 2k } $. From this, for all $s > 0$, one easily deduces a canonical isomorphism (\cite{SnTri}p.1267) of the form
\[   \begin{array}{ll}
 E_{2}^{s,t} &  \cong  Ext_{B}^{s,t}( \Sigma^{ 2l - 2k + \alpha(l) - \alpha(k) }
I^{ 2l  - 2k - \alpha(l) + \alpha(k)  }  , \mathbb{Z}/2) \\
\\
 &  \cong  
 Ext_{B}^{s + 2l  - 2k - \alpha(l) + \alpha(k) ,t  - 2l + 2k - \alpha(l) + \alpha(k) }(
   \mathbb{Z}/2 , \mathbb{Z}/2) . 
\end{array}  \]
Also there is an algebra isomorphism of the form $ Ext_{B}^{* , * }(\mathbb{Z}/2 , \mathbb{Z}/2) \cong \mathbb{Z}/2[a,b]$ where $a \in 
  Ext_{B}^{1 , 1 }$, $b \in   Ext_{B}^{1 , 3}$.
As explained in (\cite{SnTri} p.1270) $i_{k,l}$ is represented in 
\[ E_{2}^{4(k-l) + \alpha(l) - \alpha(k) , 4(k-l) + \alpha(l) - \alpha(k)}
\cong   Ext_{B}^{2k- 2l ,6k- 6l)  }(
   \mathbb{Z}/2 , \mathbb{Z}/2)   \cong \mathbb{Z}/2 = \langle b^{2k - 2l} \rangle .  \]
\end{em}
\end{dummy}
\begin{proposition}{$_{}$}
\label{3.2}

For $l < k$, in the notation of \S2, the homomorphism
\[ ( \iota_{k,l})_{*}  :   G_{4k,k}  \longrightarrow  G_{4k,l} \]
satisfies $( \iota_{k,l})_{*} (z_{4k}) = \mu_{4k,4l} 2^{2k - 2l - \alpha(k) + \alpha(l)} u^{2k-2l} z_{4l}$
for some $2$-adic unit $\mu_{4k,4l}$.
\end{proposition}
\vspace{2pt}

{\bf Proof}
\vspace{2pt}

Let $\tilde{z}_{4k} \in \tilde{G}_{4k,k}$ be as in \S\ref{2.9} so that, proved in a similar manner to 
Lemma \ref{2.10},  $2 \tilde{z}_{4k}$ is
represented in 
$E_{2}^{1,4k+1}$ in the spectral sequence
\[  \begin{array}{l}
E_{2}^{s,t} = Ext_{B}^{s,t}(  H^{*}(F_{4k}/F_{4k-1} ; \mathbb{Z}/2)  , \mathbb{Z}/2 )   \\
\\
\hspace{50pt}   \Longrightarrow
\pi_{t-s}(   bu   \wedge(F_{4k}/F_{4k-1})) \otimes  \mathbb{Z}_{2} . 
\end{array}  \]
where, from \S\ref{3.1}, we have

\[   E_{2}^{1,4k+1}  \cong  
 Ext_{B}^{1 + 2k  - \alpha(k) ,4k+1  - 2k  - \alpha(k)  }(
   \mathbb{Z}/2 , \mathbb{Z}/2) \cong \mathbb{Z}/2 = \langle a^{2k+1 - \alpha(k)} \rangle  . \]
The multiplicative pairing between these spectral sequences shows that 
\linebreak
$(\iota_{k,l})^{*}( 2 \tilde{z}_{4k}) \in \tilde{G}_{4k,l}$
is represented in the spectral
sequence
\[  \begin{array}{l}
E_{2}^{s,t} = Ext_{B}^{s,t}(  H^{*}(F_{4l}/F_{4l-1} ; \mathbb{Z}/2)  , \mathbb{Z}/2 )   \\
\\
\hspace{50pt}   \Longrightarrow
\pi_{t-s}(   bu   \wedge(F_{4l}/F_{4l-1})) \otimes  \mathbb{Z}_{2} . 
\end{array}  \]
by the generator of 
$E_{2}^{1 + 4k - 4l - \alpha(k) + \alpha(l), 1 + 8k - 4l - \alpha(k) + \alpha(l)}$
because $a^{2k+1 - \alpha(k)} b^{2k - 2l}$ is the generator of
\[   
 E_{2}^{1 + 4k - 4l - \alpha(k) + \alpha(l), 1 + 8k - 4l - \alpha(k) + \alpha(l)}  \cong  
 Ext_{B}^{1 + 4k - 2l - \alpha(k)  ,1 + 8k - 6l - \alpha(k)  }(
   \mathbb{Z}/2 , \mathbb{Z}/2)   .  \]

Since multiplication by $a$ and $b$ in the spectral sequence corresponds to multiplication by $2$ and $u$ respectively on homotopy
groups we have the following table of representatives in $\pi_{*}(   bu   \wedge(F_{4l}/F_{4l-1})) \otimes  \mathbb{Z}_{2}$.
\[ \begin{array}{c|c|c}
{\rm homotopy \ element}  &  {\rm representative}  & {\rm dimension}  \\
\hline
2z_{4l}  &  a^{2l - \alpha(l) + 1} &  4l  \\
\hline
(u/2)(2z_{4l}) &  a^{2l - \alpha(l)}b & 4l+2  \\  
\hline
(u/2)^{2}(2z_{4l}) &  a^{2l - \alpha(l) - 1}b^{2} & 4l+4  \\  
\hline
\vdots  & \vdots & \vdots  \\
\hline
(u/2)^{2l - \alpha(l)}(2z_{4l}) &  ab^{2l - \alpha(l)} & 8l - 2 \alpha(l)  \\  
\hline
u(u/2)^{2l - \alpha(l)}(2z_{4l}) &  b^{2l - \alpha(l)+1} & 8l - 2 \alpha(l) +2  \\  
\hline
u^{2}(u/2)^{2l - \alpha(l)}(2z_{4l}) &  b^{2l - \alpha(l) + 2} & 8l - 2 \alpha(l) + 4 \\  
\hline
\vdots  & \vdots & \vdots  \\
\hline
\end{array} \]

Therefore there are two cases for $(\iota_{k,l})_{*}(2 \tilde{z}_{4k})$.
If $2k - 2l \geq 2l - \alpha(l) + 1$ then $b^{2k-2l}$ represents
$ u^{2k-2l - (2l - \alpha(l))} (u/2)^{2l - \alpha(l)} \tilde{z}_{4l} = 
 u^{2k-4l + \alpha(l))} (u/2)^{2l - \alpha(l)} \tilde{z}_{4l}  $
and, up to multiplication by $2$-adic units,
$ (\iota_{k,l})_{*}(2 \tilde{z}_{4k})$ is equal to 
\linebreak
$ 2^{1 + 2k - \alpha(k)} u^{2k-4l + \alpha(l))} (u/2)^{2l - \alpha(l)} \tilde{z}_{4l}$, as required.   
On the other hand, if 
\linebreak
$2k - 2l \leq 2l - \alpha(l) $ then
$a^{2l - \alpha(l) + 1 - (2k-2l)}b^{2k-2l} = a^{4l - 2k - \alpha(l) + 1}b^{2k-2l}$ represents
$(u/2)^{2k-2l}(2 \tilde{z}_{4l}) $ which shows that, up to $2$-adic units, 
$ (\iota_{k,l})_{*}(2 \tilde{z}_{4k}) $ is equal to
$ 2^{1 + 2k - \alpha(k) - (4l - 2k - \alpha(l) + 1)} (u/2)^{2k- 2l }(2 \tilde{z}_{4l})
=   2^{4k - \alpha(k) - 4l  + \alpha(l) } (u/2)^{2k- 2l }(2 \tilde{z}_{4l}) $, as required. $\Box$ 
\begin{proposition}{$_{}$}
\label{3.3}

Let $\psi^{3} : bo \longrightarrow  bo$ denote the Adams operation, as usual. Then, in the notation of \S\ref{2.7},
\[ (1 \wedge \psi^{3})_{*}(g_{4k,4k}) = \left\{
 \begin{array}{ll}
 9^{k} g_{4k,4k} +  9^{k-1}  2^{  \nu_{2}(k)  + 3  } g_{4k, 4k-4} & {\rm if}  \  k \geq 3,  \\
\\
 9^{2} g_{8,8} +  9 \cdot 2^{ 3  } g_{8, 4} & {\rm if}  \  k = 2, \\
\\
 9 g_{4,4} +  2 g_{4, 0} & {\rm if}  \  k = 1,  \\
\\
g_{0,0}  & {\rm if}  \  k = 0.  
\end{array}  \right.  \]
\end{proposition}
\vspace{2pt}

\newpage
{\bf Proof}
\vspace{2pt}

The map $ (1 \wedge \psi^{3})_{*}$ fixes $u$, multiplies $v$ by $9$ and is multiplicative. Therefore
\[  \begin{array}{ll}

 (1 \wedge \psi^{3})_{*}(c_{4k}) &
 =  \Pi_{i=1}^k  \big(\frac{9 v^2-9^{i-1}u^2}{9^k-9^{i-1}}\big)  \\
\\&  
=  9^{k-1}   \big(\frac{ (9v^{2} - 9^{k}u^2 + 9^{k}u^2  - u^2) \Pi_{i=2}^k ( v^2-9^{i-2}u^2)}{   \Pi_{i=1}^k (9^k-9^{i-1})}\big) \\
\\
& =  9^{k-1}   \big(\frac{ (9v^{2} - 9^{k}u^2 ) \Pi_{i=2}^k ( v^2-9^{i-2}u^2)}{   \Pi_{i=1}^k (9^k-9^{i-1})}\big)
+  9^{k-1}   \big(\frac{ ( 9^{k}u^2  - u^2) \Pi_{i=2}^k ( v^2-9^{i-2}u^2)}{   \Pi_{i=1}^k (9^k-9^{i-1})}\big)  \\
\\
& =   9^{k}   \big(\frac{ (v^{2} - 9^{k-1}u^2 ) \Pi_{i=2}^k ( v^2-9^{i-2}u^2)}{   \Pi_{i=1}^k (9^k-9^{i-1})}\big)
+  9^{k-1}   \big(\frac{ ( 9^{k}u^2  - u^2) \Pi_{i=2}^k ( v^2-9^{i-2}u^2)}{   \Pi_{i=1}^k (9^k-9^{i-1})}\big)  \\
\\
& =   9^{k} c_{4k}
+  9^{k-1} (9^{k}-1)  \big(\frac{  u^2 \Pi_{i=1}^{k-1} ( v^2-9^{i-1}u^2)}{  (9^{k}-1)  \Pi_{i=1}^{k-1} (9^k-9^{i-1})}\big)  \\
\\
& =   9^{k} c_{4k}
+  9^{k-1}  u^2  c_{4k-4}.
\end{array} \]

Hence, for $k \geq 1$, we have
\[  \begin{array}{ll}
 (1 \wedge \psi^{3})_{*}(f_{4k}) &
 =  2^{2k - \alpha(k) }  (1 \wedge \psi^{3})_{*}(c_{4k})  \\
\\
& =  2^{2k - \alpha(k) }  9^{k} c_{4k}
+  9^{k-1} u^2  2^{2k - \alpha(k) - 2k + 2 + \alpha(k-1) + 2k-2 - \alpha(k-1)} c_{4k-4}  \\
\\
& =  9^{k} f_{4k} +  9^{k-1} u^2  2^{2 - \alpha(k)  + \alpha(k-1) } f_{4k-4} \\
\\
& =  9^{k} f_{4k} +  9^{k-1} u^2  2^{  \nu_{2}(k)  + 1  } f_{4k-4} ,
\end{array} \]
which yields the result, by the formulae of \S\ref{2.7}.  $\Box$
\begin{proposition}{$_{}$}
\label{3.4}

When $k > l$ 
\[ (1 \wedge \psi^{3})_{*}(g_{4k,4l}) = \left\{
 \begin{array}{ll}
  9^{l} g_{4k,4l} + 9^{l-1}  g_{4k,4l-4}  & {\rm if}  \  4l - \alpha(l) \leq 2k,  \\
\\
 9^{l} g_{4k,4l}  +   9^{l-1} 2^{4l - \alpha(l) - 2k}  g_{4k,4l-4}   &  {\rm if }  \  4l - \alpha(l) - \nu_{2}(l) - 3   \\
&  \hspace{40pt} \leq 2k <  4l - \alpha(l),  \\
\\
  9^{l} g_{4k,4l}  +   9^{l-1} 2^{3 + \nu_{2}(k)}  g_{4k,4l-4}   &   {\rm if }  \  2k < 4l - \alpha(l) - \nu_{2}(l) - 3 \\
&  \hspace{40pt} < 4l - \alpha(l)
\end{array}  \right.  \]
\end{proposition}
\vspace{2pt}

\newpage
{\bf Proof}
\vspace{2pt}

Suppose that $4l - \alpha(l) \leq 2k$ then, by Proposition \ref{3.3} (proof),
\[ \begin{array}{l}
 (1 \wedge \psi^{3})_{*}(g_{4k,4l})  \\
\\
=   (1 \wedge \psi^{3})_{*}(  u^{2k-4l+\alpha(l)}\big[\frac{u^{2l-\alpha(l)}f_{4l}}{2^{2l-\alpha(l)}}\big]  )  \\
\\
=   u^{2k-4l+\alpha(l)}\big[\frac{u^{2l-\alpha(l)}(  9^{l} f_{4l} +  9^{l-1} u^2  2^{  \nu_{2}(l)  + 1  }
f_{4l-4})}{2^{2l-\alpha(l)}}\big]   \\
\\
=   9^{l} g_{4k,4l} +   9^{l-1} u^{2k-4l+\alpha(l)}\big[\frac{u^{2l-\alpha(l)}   u^2  2^{  \nu_{2}(l)  + 1  }
f_{4l-4}}{2^{2l-\alpha(l)}}\big]   \\
\\
=   9^{l} g_{4k,4l} +   9^{l-1} u^{2k-4l+\alpha(l)}\big[\frac{u^{2l+2-\alpha(l)}  2^{  \nu_{2}(l)  + 1  }
f_{4l-4}}{2^{2l-\alpha(l)}}\big] .
\end{array}  \]
Then, since $\nu_{2}(l) = 1 + \alpha(l-1) - \alpha(l)$, 
\[ 4(l-1) - \alpha(l-1) = 4l - \alpha(l) + \alpha(l) - \alpha(l-1) -4  =  4l - \alpha(l)   - 3 -
\nu_{2}(l)  < 2k \]
so that
\[ \begin{array}{ll}
 g_{4k,4l-4} &  =    
u^{2k-4l + 4 +\alpha(l-1)}\big[\frac{u^{2l-2-\alpha(l-1)}f_{4l-4}}{2^{2l-2-\alpha(l-1)}}\big]  \\
\\
& = u^{2k-4l+\alpha(l)}\big[\frac{u^{2l+2 -\alpha(l)}f_{4l-4}}{2^{2l-2 - \alpha(l) + \alpha(l)-\alpha(l-1)}}\big]  \\
\\
& = u^{2k-4l+\alpha(l)}\big[\frac{u^{2l+2 -\alpha(l)}f_{4l-4}}{2^{2l - \alpha(l) - \nu_{2}(l)  - 1}}\big] 
\end{array}   \]
so that, for $0 < l < k$ suppose that $4l - \alpha(l) \leq 2k$,
\[  (1 \wedge \psi^{3})_{*}(g_{4k,4l}) = 9^{l} g_{4k,4l} + 9^{l-1}  g_{4k,4l-4}  .  \]

Similarly, for $0 < l < k$ if $4l - \alpha(l) > 2k$ then, by Proposition \ref{3.3} (proof),
\[ \begin{array}{l}
 (1 \wedge \psi^{3})_{*}(g_{4k,4l})  \\
\\
=   (1 \wedge \psi^{3})_{*}( \big[\frac{u^{2(k-l)f_{4l}}}{2^{2(k-l)}}\big] )  \\
\\
=  \big[\frac{u^{2(k-l) (  9^{l} f_{4l} +  9^{l-1} u^2  2^{  \nu_{2}(k)  + 1  } f_{4k-4} ) }}{2^{2(k-l)}}\big]  \\
\\
=  9^{l} g_{4k,4l}  +   9^{l-1}  \big[ \frac{ u^{2k-2l+2}   2^{  \nu_{2}(k)  + 1  } f_{4k-4}  }{2^{2(k-l)} }\big].
\end{array} \]
This situation splits into two cases given by

(i)  \  $4l - \alpha(l) - \nu_{2}(l) - 3 \leq 2k <  4l - \alpha(l)$ or

(ii)  \  $2k < 4l - \alpha(l) - \nu_{2}(l) - 3  < 4l - \alpha(l) $.

In case (i) $4l - 4 - \alpha(l-1) = 4l - \alpha(l) - \nu_{2}(l) - 3 \leq 2k$ and so again we have
\[ \begin{array}{ll}
 g_{4k,4l-4} &  =    
u^{2k-4l + 4 +\alpha(l-1)}\big[\frac{u^{2l-2-\alpha(l-1)}f_{4l-4}}{2^{2l-2-\alpha(l-1)}}\big]  \\
\\
& = \frac{u^{2k - 2l + 2}f_{4l-4}}{2^{2l-1 - \nu_{2}(l) - \alpha(l) }}  \\
\\
& = \frac{u^{2k - 2l + 2}2^{1 + \nu_{2}(l)}f_{4l-4}}{2^{2k - 2l + 4l - \alpha(l) - 2k}}
\end{array}   \]
so that
\[ (1 \wedge \psi^{3})_{*}(g_{4k,4l})  =   9^{l} g_{4k,4l}  +   9^{l-1} 2^{4l - \alpha(l) - 2k}  g_{4k,4l-4} .\]

In case (ii)
\[ \begin{array}{ll}
 g_{4k,4l-4} &  =   \big[ \frac{ u^{2k-2l+2}   f_{4k-4}  }{2^{2(k-l+2)} }\big] 
\end{array}  \]
so that
\[ (1 \wedge \psi^{3})_{*}(g_{4k,4l})  =   9^{l} g_{4k,4l}  +   9^{l-1} 2^{3 + \nu_{2}(k)}  g_{4k,4l-4} .\]
$\Box$
\begin{dummy}
\label{3.5}
\begin{em}

In the notation of \S\ref{2.1}, suppose that $A \in  U_{\infty}\mathbb{Z}_2$ satisfies
\[  \Psi(A^{-1})  = [1 \wedge \psi^{3}]  \in  {\rm Aut}_{left-bu-mod}^{0}( bu \wedge bo )    .  \]

Therefore, by definition of $\Psi$ and the formula of Theorem \ref{2.11}
\[  \begin{array}{ll}
\sum_{l \leq k} A_{l,k}  (\iota_{k,l})_{*}(z_{4k}) & =  (1 \wedge \psi^{3})_{*}(z_{4k}) \\
\\
& =   \Sigma_{i=0}^k2^{\beta(k,i)}  \lambda_{4k,4i} (1 \wedge \psi^{3})_{*}(g_{4k,4i} ).
\end{array} \]
On the other hand
\[ \begin{array}{l}
\sum_{l \leq k} A_{l,k}  (\iota_{k,l})_{*}(z_{4k}) \\
\\
= A_{k,k} z_{4k} +  \sum_{l < k}  A_{l,k} \mu_{4k,4l} 2^{2k - 2l - \alpha(k) + \alpha(l)} u^{2k-2l} z_{4l} \\
\\
= A_{k,k}  \Sigma_{i=0}^k2^{\beta(k,i)}  \lambda_{4k,4i} g_{4k,4i}  \\
\\
\hspace{30pt} +  \sum_{l < k} \  \Sigma_{i=0}^{l} \  A_{l,k} \  \mu_{4k,4l} 2^{2k - 2l - \alpha(k) + \alpha(l)} u^{2k-2l} 
 2^{\beta(l,i)}  \lambda_{4l,4i}  g_{4l,4i} .
\end{array} \]

In order to determine the $ A_{k,l} $'s it will suffice to express
$u^{2k-2l} g_{4l,4i}$ as a multiple of $g_{4k,4i}$ and then to equate coefficients in the above expressions.
By definition
\[ \begin{array}{ll}
u^{2k-2l} g_{4l,4i} & =    \left\{
\begin{array}{ll}
u^{2k-2l}u^{2l-4i+\alpha(i)}\big[\frac{u^{2i-\alpha(i)}f_{4i}}{2^{2i-\alpha(i)}}\big] &   {\rm if}  \   4i-\alpha(i) \leq 2l ,  \\
\\
u^{2k-2l} \big[\frac{u^{2(l-i)f_{4i}}}{2^{2(l-i)}}\big]  &  {\rm if}  \   4i-\alpha(i) > 2l   . 
\end{array}  \right. \\
\\
&  =    \left\{
\begin{array}{ll}
 \frac{u^{2k- 2i}f_{4i}}{2^{2i-\alpha(i)}} &   {\rm if}  \   4i-\alpha(i) \leq 2l ,  \\
\\
 \frac{u^{2k-2i}f_{4i}}{2^{2l- 2i}}  &  {\rm if}  \   4i-\alpha(i) > 2l   
\end{array} \right.
\end{array} \]
while
\[ g_{4k,4i} =    \left\{
\begin{array}{ll}
u^{2k-4i+\alpha(i)}\big[\frac{u^{2i-\alpha(i)}f_{4i}}{2^{2i-\alpha(i)}}\big] &   {\rm if}  \   4i-\alpha(i) \leq 2k ,  \\
\\
\big[\frac{u^{2(k-i)f_{4i}}}{2^{2(k-i)}}\big]  &  {\rm if}  \   4i-\alpha(i) > 2k   .
\end{array} \right.   \]
From these formulae we find that
\[ u^{2k - 2l} g_{4l,4i} = \left\{
\begin{array}{ll}
g_{4k,4i}  & {\rm if}  \  4i - \alpha(i) \leq 2l \leq  2k ,  \\
\\
  2^{4i - \alpha(i) - 2l} g_{4k,4i} &  {\rm if}  \  2l < 4i - \alpha(i)  \leq  2k ,  \\
\\
 2^{2k - 2l} g_{4k,4i}  &  {\rm if}  \  2l < 2k < 4i - \alpha(i) . 
\end{array}  \right.  \]

Now let us calculate $A_{l,k}$.

When $k=0$ we have $z_{0} = (1 \wedge \psi^{3})_{*}(z_{0}) = A_{0,0} (\iota_{0,0})_{*}(z_{0}) = A_{0,0} z_{0}$
so that $A_{0,0} = 1$.

When $k=1$ we have
\[  \begin{array}{ll}
\sum_{l \leq 1} A_{l,1}  (\iota_{1,l})_{*}(z_{4}) & =  (1 \wedge \psi^{3})_{*}(z_{4}) \\
\\
& =    \lambda_{4,4} (1 \wedge \psi^{3})_{*}(g_{4,4} ) + 2 \lambda_{4,0}  (1 \wedge \psi^{3})_{*}(g_{4,0} )  \\
\\
& =  \lambda_{4,4} ( 9g_{4,4}  + 2  g_{4,0}) +  2 \lambda_{4,0}g_{4,0} 
\end{array} \]
and
\[  \begin{array}{ll}
\sum_{l \leq 1} A_{l,k}  (\iota_{1,l})_{*}(z_{4}) 
& =  A_{1,1} z_{4} + A_{0,1} \mu_{1,0} 2 g_{4,0}   \\
\\
& =   A_{1,1} (2 \lambda_{4,0} g_{4,0} + \lambda_{4,4} g_{4,4}) + A_{0,1} \mu_{1,0} 2 g_{4,0} 
\end{array} \]
which implies that $A_{1,1} = 9$ and
$ A_{0,1} = \mu_{1,0}^{-1} ( \lambda_{4,4}    -  8 \lambda_{4,0})  $  so that $A_{0,1} \in \mathbb{Z}_{2}^{*}$.

When $k=2$ we have
\[  \begin{array}{ll}
\sum_{l \leq 2} A_{l,2}  (\iota_{2,l})_{*}(z_{8}) & =  (1 \wedge \psi^{3})_{*}(z_{8}) \\
\\
& =    (1 \wedge \psi^{3})_{*}(  \lambda_{8,8} g_{8,8} +  2^{3}  \lambda_{8,4} g_{8,4} +  2^{3}  \lambda_{8,0} g_{8,0})   \\
\\
& =     \lambda_{8,8} (9^{2}g_{8,8} + 9 \cdot 2^{3} g_{8,4}) +  2^{3}  \lambda_{8,4}(9 g_{8,4} + g_{8,0}) +  2^{3}  \lambda_{8,0}
g_{8,0} 
\end{array} \]
and
\[  \begin{array}{ll}
\sum_{l \leq 2} A_{l,2}  (\iota_{2,l})_{*}(z_{8}) & = A_{2,2} z_{8} +  A_{1,2}  (\iota_{2,1})_{*}(z_{8}) 
+ A_{0,2}  (\iota_{2,0})_{*}(z_{8})  \\
\\
& = A_{2,2}(   \lambda_{8,8} g_{8,8} +  2^{3}  \lambda_{8,4} g_{8,4} +  2^{3}  \lambda_{8,0} g_{8,0} )  \\
\\
&  \hspace{20pt} +   A_{1,2}  ( \mu_{8,4} 2^{2} u^{2} z_{4}) 
+ A_{0,2}  (\mu_{8,0} 2^{3} u^{4}z_{0})    \\
\\
& =   A_{2,2}(   \lambda_{8,8} g_{8,8} +  2^{3}  \lambda_{8,4} g_{8,4} +  2^{3}  \lambda_{8,0} g_{8,0} )  \\
\\
& \hspace{20pt} +   A_{1,2}  \mu_{8,4} 2^{2}(2 \lambda_{4,0} g_{8,0} + \lambda_{4,4}  u^{2}  g_{4,4}) 
+ A_{0,2} \mu_{8,0} 2^{3} g_{8,0} \\
\\
& =   A_{2,2}(   \lambda_{8,8} g_{8,8} +  2^{3}  \lambda_{8,4} g_{8,4} +  2^{3}  \lambda_{8,0} g_{8,0} )  \\
\\
& \hspace{20pt} +   A_{1,2}  \mu_{8,4} 2^{2}(2 \lambda_{4,0} g_{8,0} + \lambda_{4,4}  2  g_{8,4}) 
+ A_{0,2} \mu_{8,0} 2^{3} g_{8,0} .
\end{array} \]
Therefore we obtain
\[ \begin{array}{l}
  \lambda_{8,8} (9^{2}g_{8,8} + 9 \cdot 2^{3} g_{8,4}) +  2^{3}  \lambda_{8,4}(9 g_{8,4} + g_{8,0}) +  2^{3} 
\lambda_{8,0} g_{8,0} \\
\\
 =  A_{2,2}(   \lambda_{8,8} g_{8,8} +  2^{3}  \lambda_{8,4} g_{8,4} +  2^{3}  \lambda_{8,0} g_{8,0} )  \\
\\
\hspace{40pt} +   A_{1,2}  \mu_{8,4}
2^{2}(2 \lambda_{4,0} g_{8,0} + \lambda_{4,4}  2  g_{8,4})  + A_{0,2} \mu_{8,0} 2^{3} g_{8,0}
\end{array}  \]
which yields
\[ \begin{array}{l}
   9^{2}    =  A_{2,2}  ,   \\
\\
 \lambda_{8,8} \cdot  9   +  \lambda_{8,4}( 9 -  9^{2})      =   A_{1,2} \mu_{8,4} \lambda_{4,4}   , \\
\\
     \lambda_{8,4}  +   \lambda_{8,0} (1 - 9^{2}) =    
  A_{1,2}  \mu_{8,4} \lambda_{4,0}   + A_{0,2} \mu_{8,0}  .
\end{array}  \]
Hence $A_{1,2} \in \mathbb{Z}_{2}^{*}$. 

Now assume that $k \geq 3$ and consider the relation derived above
\[ \begin{array}{l}
   \Sigma_{i=0}^k  2^{\beta(k,i)}  \lambda_{4k,4i} (1 \wedge \psi^{3})_{*}(g_{4k,4i} ) \\
\\
= A_{k,k}  \Sigma_{i=0}^k 2^{\beta(k,i)}  \lambda_{4k,4i} g_{4k,4i}  \\
\\
\hspace{30pt} +  \sum_{l < k} \  \Sigma_{i=0}^{l} \  A_{l,k} \  \mu_{4k,4l} 2^{2k - 2l - \alpha(k) + \alpha(l)} u^{2k-2l} 
 2^{\beta(l,i)}  \lambda_{4l,4i}  g_{4l,4i} .
\end{array} \]
The coefficient of $g_{4k,4k}$ on the left side of this relation is equal to
$ \lambda_{4k,4k} 9^{k}$ and on the right side it is $A_{k,k} \lambda_{4k,4k}$ so that $A_{k,k} = 9^{k}$ for all $k \geq 3$.
From the coefficient of $g_{4k,4k-4}$ we obtain the relation
\[  \begin{array}{l}
 \lambda_{4k,4k}  9^{k-1}  2^{  \nu_{2}(k)  + 3  } +
2^{3 + \nu_{2}(k)}  \lambda_{4k,4k-4} 9^{k-1} \\
\\
 = 9^{k} 2^{3 + \nu_{2}(k)}  \lambda_{4k,4k-4}  \\
\\
\hspace{30pt}   
+  A_{k-1,k} \  \mu_{4k,4k-4} 2^{2 - \alpha(k) + \alpha(k-1)} 2^{2} \lambda_{4k-4,4k-4}  
2^{3 + \nu_{2}(k)}  \lambda_{4k,4k-4} 9^{k-1} \\
\\
 = 9^{k} 2^{3 + \nu_{2}(k)}  \lambda_{4k,4k-4}  \\
\\
\hspace{30pt}   
+  A_{k-1,k} \  \mu_{4k,4k-4} 2^{3 + \nu_{2}(k)} \lambda_{4k-4,4k-4}  
\end{array} \]
which shows that $A_{k-1,k} \in \mathbb{Z}_{2}^{*}$ for all $k \geq 3$.
This means that we may conjugate $A$ by the matrix
$D = {\rm diag}(1 , A_{1,2}, A_{1,2}A_{2,3} , A_{1,2}A_{2,3} A_{3,4}, \ldots ) \in U_{\infty}\mathbb{Z}_{2}$ to obtain
\[   D A D^{-1} =  C = \left( 
 \begin{array}{cccccc} 

1 & 1 & c_{1,3} & c_{1,4} & c_{1,5} & \ldots \\
\\
0 & 9 & 1 & c_{2,4} & c_{2,5} & \ldots \\
\\
0 & 0 & 9^{2} & 1 & c_{3,5}  & \ldots \\
\\
0 & 0 & 0 & 9^{3} & 1   & \ldots \\
\\
\vdots &  \vdots &  \vdots &  \vdots &  \vdots &  \vdots 
\end{array}  \right)  . \]

In the next section we examine whether we can conjugate this matrix further in $U_{\infty}\mathbb{Z}_{2}$
to obtain the matrix
\[  B = \left( 
 \begin{array}{cccccc} 

1 & 1 & 0 & 0 & 0 & \ldots \\
\\
0 & 9 & 1 & 0 & 0 & \ldots \\
\\
0 & 0 & 9^{2} & 1 & 0  & \ldots \\
\\
0 & 0 & 0 & 9^{3} & 1   & \ldots \\
\\
\vdots &  \vdots &  \vdots &  \vdots &  \vdots &  \vdots 
\end{array}  \right).\]
\end{em}
\end{dummy}

\section{The Matrix Reloaded}
\begin{dummy}
\label{4.1}
\begin{em}

Let $B , C \in U_{\infty}\mathbb{Z}_{2}$ denote the upper triangular matrices which occurred in \S\ref{3.5}
\[  B = \left( 
 \begin{array}{cccccc} 

1 & 1 & 0 & 0 & 0 & \ldots \\
\\
0 & 9 & 1 & 0 & 0 & \ldots \\
\\
0 & 0 & 9^{2} & 1 & 0  & \ldots \\
\\
0 & 0 & 0 & 9^{3} & 1   & \ldots \\
\\
\vdots &  \vdots &  \vdots &  \vdots &  \vdots &  \vdots 
\end{array}  \right) , 
C = \left( 
Ê\begin{array}{cccccc} 

1 & 1 & c_{1,3} & c_{1,4} & c_{1,5} & \ldots \\
\\
0 & 9 & 1 & c_{2,4} & c_{2,5} & \ldots \\
\\
0 & 0 & 9^{2} & 1 & c_{3,5} Ê& \ldots \\
\\
0 & 0 & 0 & 9^{3} & 1 Ê & \ldots \\
\\
\vdots & Ê\vdots & Ê\vdots & Ê\vdots & Ê\vdots & Ê\vdots 
\end{array} Ê\right) . \]

The following result is the main result of this section. Along with the discussion of \S\ref{3.5} it completes the proof of Theorem \ref{0.1}.
\begin{theorem}{$_{}$}
\label{4.2}
There exists an upper triangular matrix $U  \in U_{\infty}\mathbb{Z}_{2}$ such that $U^{-1}CU = B$.
\end{theorem}
\vspace{2pt}

{\bf Proof}
\vspace{2pt}

Let $U$ have the form
\[  U = \left(
\begin{array}{ccccc}
1 & u_{1,2} & u_{1,3} & u_{1,4} & \ldots \\
\\
0 & 1 + (9-1)u_{1,2} & u_{2,3} & u_{2,4} & \ldots \\
\\
0 & 0 & 1 + (9-1)u_{1,2} + (9^{2}-9)u_{2,3} & Ê u_{3,4} & \ldots \\
\\
\ldots & Ê\ldots & Ê\ldots & Ê\ldots & Ê\ldots Ê \\
\\
\ldots & Ê\ldots & Ê\ldots & Ê\ldots & Ê\ldots
\end{array} \right) Ê .  \]
Then $(UB)_{j,j} = U_{j,j} B_{j,j} = Ê C_{j,j} U_{j,j} = (CU)_{j,j}$ and, in fact, 
$(UB)_{j,j+1} = (CU)_{j,j+1}$ for all $j$, too. ÊFor any $1<s<j$
we have 
\[ \begin{array}{l}
 (UB)_{j-s,j}=u_{j-s,j}9^{j-1}+u_{j-s,j-1}  \  {\rm and}  \\
\\ 
 (CU)_{j-s,j}=9^{j-s-1}u_{j-s,j}+u_{j-s+1,j}+c_{j-s,j-s+2}u_{j-s+2,j}+\cdots+c_{j-s,j}u_{j,j}  . 
\end{array} \]
In order to prove Theorem \ref{4.2} it suffices to verify that we are able to solve for the $u_{i,j}$ in
the equations  $(UB)_{s,t}=(CU)_{s,t}$ for all
$s\leq t$ inductively in such in a manner such that, for every $k$, the first $k$-columns of the equality $UB = CU$ 
is achieved after a finite number of steps.
Lemma \ref{4.3} provides a method which proceeds inductively on the columns of $U$. $\Box$

\begin{lemma}{$_{}$}
\label{4.3}
For $j\geq3$ and $1<s<j$, $u_{j-s,j-1}$ may be written as a linear combination of $u_{j-2,j-1},u_{j-3,j-1},\ldots,u_{j-s+1,j-1}$ and $u_{j-1,j},u_{j-2,j},\ldots,u_{1,j}$
\end{lemma}
\vspace{2pt}

{\bf Proof}
\vspace{2pt}

We shall prove the result by induction on $j$. ÊConsider the case $j=3$, we have the following equation:
$$u_{3-s,3}9^2+u_{3-s,2}=9^{2-s}u_{3-s,3}+u_{4-s,3}+c_{3-s,5-s}u_{5-s,3}+\cdots+c_{3-s,3}u_{3,3}$$ for $1<s<3\implies s=2$. Ê
Hence substituting $s=2$ gives\\

$u_{1,3}9^2+u_{1,2}=u_{1,3}+u_{2,3}+c_{1,3}u_{3,3}$\\

$\implies u_{1,2}=(1-9^2)u_{1,3}+u_{2,3}+c_{1,3}(1+(9-1)u_{1,2}+(9^2-9)u_{2,3})$\\

$\implies (1-(9-1)c_{1,3})u_{1,2}=(1-9^2)u_{1,3}+(1+(9^2-9))c_{1,3}u_{2,3}+c_{1,3}$\\
and since $(1-(9-1)c_{1,3})$ is a 2-adic unit we can write $u_{1,2}$ as a $\mathbb{Z}_{2}$-linear combination of $u_{1,3}$ and $u_{2,3}$
as required.\\

We now need to show that if the lemma is true for Ê$j=3,4,\ldots,k-1$ then it is also true for $j=k$. ÊThis means we need to solve
\[  \begin{array}{l}
u_{k-s,k}9^{k-1}+u_{k-s,k-1}  \\
\\
=9^{k-s-1}u_{k-s,k}+u_{k-s+1,k}+c_{k-s,k-s+2}u_{k-s+2,k}+\cdots+c_{k-s,k}u_{k,k} 
 \end{array} \] 
for $u_{k-s,k-1}$ for $1<s<k$.
ÊThis equation may be rewritten
\[  \begin{array}{l}
u_{k-s,k-1}  \\
\\
 =    (9^{k-s-1}-9^{k-1})u_{k-s,k}+u_{k-s+1,k}+c_{k-s,k-s+2}u_{k-s+2,k}+\cdots+c_{k-s,k}u_{k,k}  \\
\\
=  (9^{k-s-1}-9^{k-1})u_{k-s,k}+u_{k-s+1,k}+c_{k-s,k-s+2}u_{k-s+2,k}+\cdots  \\
\\
\cdots  +  c_{k-s,k}(1+(9-1)u_{1,2}+(9^2-9)u_{2,3}+\cdots+(9^{k-1}-9^{k-2})u_{k-1,k}).
\end{array}  \]
Now consider the case $s=k-1$
\[  \begin{array}{l}
 u_{1,k-1}  \\
\\
=  (1-9^{k-1})u_{1,k}+u_{2,k}+c_{1,3}u_{3,k} + \cdots \\
\\
\cdots  c_{1,k}\underbrace{(1+(9-1)u_{1,2}+(9^2-9)u_{2,3}+\cdots+(9^{k-1} - 9^{k-2})u_{k-1,k})}_{B}
\end{array}   \]
By repeated substitutions the bracket $B$ may be rewritten as a linear combination of 
$u_{1,k-1},u_{2,k-1},\ldots,u_{k-2,k-1}$ and $u_{k-1,k}$. ÊThe important point to notice about this linear combination is that the
coefficient of
$u_{1,k-1}$ will be an even 2-adic integer. Hence, we can move this term to the left hand side of the equation to obtain a 2-adic unit
times $u_{1,k-1}$ equals a linear combination of $u_{1,k},u_{2,k},\ldots,u_{k-1,k}$ and $u_{2,k-1},u_{3,k-1},\ldots,u_{k-2,k-1}$ as
required.

Now consider $s=k-2$
\[  \begin{array}{l}
 u_{2,k-1} \\
\\
 =  (9-9^{k-1})u_{2,k}+u_{3,k}+c_{2,4}u_{4,k}+  \cdots  \\
\\
  \cdots  +  c_{2,k}\underbrace{(1+(9-1)u_{1,2}+(9^2-9)u_{2,3}+\cdots+(9^{k-1}-9^{k-2})u_{k-1,k}).}_{B'}
\end{array}  \]

As before $B'$ can be written as a linear combination of 
\linebreak
$u_{1,k-1},u_{2,k-1},\ldots,u_{k-2,k-1}, u_{k-1,k}$ and from the case
$s=k-1$,
$u_{1,k-1}$ may be replaced by a linear combination of 
$u_{1,k},\ldots,u_{k-1,k}$ and $u_{2,k-1},\ldots,u_{k-2,k-1}$. ÊAgain the
important observation is that the coefficient of $u_{2,k-1}$ is an even 2-adic integer, hence this term can be moved to the left
hand side of the equation to yield a 2-adic unit times $u_{2,k-1}$ equals a linear combination of $u_{1,k},u_{2,k},\ldots,u_{k-1,k}$
and $u_{3,k-1},\ldots,u_{k-2,k-1}$ as required.

Clearly this process may be repeated for $s=k-3,k-4,\ldots,2$ to get a 2-adic unit times $u_{k-s,k-1}$ as a linear 
combination of $u_{1,k},u_{2,k},\ldots,u_{k-1,k}$ and $u_{k-s+1,k-1},\ldots,u_{k-2,k-1}$ as required. $\Box$

\section{Applications }
\begin{dummy}
\label{5.0a}{$bu \wedge bu$}
\begin{em}

Theorem \ref{0.1} implies that, in the $2$-local stable homotopy category there exists an equivalence 
$C' \in  {\rm Aut}_{left-bu-mod}^{0}( bu \wedge bo ) $ such that
\[ C' (1 \wedge \psi^{3}) C'^{-1} = \sum_{k \geq 0} \ 9^{k} \iota_{k,k} + \sum_{k \geq 1} \ \iota_{k,k-1}    \] 
where $\iota_{k,l}$ is as in \S\ref{2.1}, considered as left $bu$-endomorphism of $bu \wedge bo$ via the equivalence $\hat{L}$ of \S\ref{2.1}.

In \cite{SnTri} use is made of an equivalence of the form $bu \simeq bo \wedge \Sigma^{-2} \mathbb{CP}^{2}$, first noticed by Reg Wood (as
remarked in \cite{A}) and independently by Don Anderson (both unpublished). This is easy to construct.
By definition $ bu^{0}(  \Sigma^{-2}  \mathbb{CP}^{2} ) \cong bu^{2}(\mathbb{CP}^{2}) \cong [ \mathbb{CP}^{2} , BU ] $
and from the cofibration sequence $S^{0} \longrightarrow   \Sigma^{-2}  \mathbb{CP}^{2}  \longrightarrow  S^{2}$ we see that
$bu^{0}( \Sigma^{-2}  \mathbb{CP}^{2} ) \cong \mathbb{Z} \oplus \mathbb{Z}$ fitting into the following exact sequence
\[  0 \longrightarrow  bu^{0}(S^{2}) \longrightarrow   bu^{0}(  \Sigma^{-2}  \mathbb{CP}^{2} ) \longrightarrow   bu^{0}( S^{0} ) \longrightarrow
0 .  \]
Choosing any stable homotopy class $x: \Sigma^{-2}  \mathbb{CP}^{2} \longrightarrow bu$ restricting to the generator of $ bu^{0}( S^{0} )$
yields an equivalence of the form 
\[  bo \wedge ( \Sigma^{-2}  \mathbb{CP}^{2}) \stackrel{c \wedge x}{\longrightarrow}  bu \wedge bu 
\stackrel{\mu}{\longrightarrow}  bu  \]
in which $c$ denoted complexification and $\mu$ is the product.

In the $2$-local stable homootopy category there is a map
\[ \Psi : \Sigma^{-2}  \mathbb{CP}^{2} \longrightarrow  \Sigma^{-2}  \mathbb{CP}^{2}  \] 
which satisfies $\Psi^{*}( z ) = \psi^{3}(z)$ for all $z \in bu^{0}(\Sigma^{-2}  \mathbb{CP}^{2})$. For example, take
$\Psi$ to be $3^{-1}$ times the double
desuspension  of the restriction to the four-skeleton of the CW complex $ \mathbb{CP}^{\infty} = BS^{1}$ of the map induced by $z \mapsto z^{3}$ on $S^{1}$, the
circle.
With this definition there is a homotopy commutative diagram in the $2$-local stable homotopy category
\newline
\begin{picture}(259,182)
\thinlines    \put(54,16){$ bo  $}
              \put(263,18){$ bu$}
              \put(245,166){$bo  \wedge \Sigma^{-2}  \mathbb{CP}^{2}$}
              \put(40,166){$bo   \wedge \Sigma^{-2}  \mathbb{CP}^{2}$}
              \put(60,143){\vector(0,-1){95}}
              \put(264,143){\vector(0,-1){92}}
              \put(120,20){\vector(1,0){88}}
              \put(120,168){\vector(1,0){88}}
              \put(160,155){$ \psi^{3} \wedge \Psi$}
              \put(160,25){$ \psi^{3}$}
              \put(65,95){$ \simeq$}
              \put(255,95){$ \simeq$}
\end{picture}
\newline
in which the vertical maps are equal given by the Anderson-Wood equivalence. 

Now suppose that we form the smash product with $ \Sigma^{-2}  \mathbb{CP}^{2}$ of the $2$-local left $bu$-module equivalence
$bu \wedge bo \simeq \vee_{k \geq 0}  bu \wedge (F_{4k}/F_{4k-1})$ to obtain a left $bu$-module equivalence of the form
\[ bu \wedge bu  \simeq \vee_{k \geq 0}  bu \wedge (F_{4k}/F_{4k-1}) \wedge   \Sigma^{-2}  \mathbb{CP}^{2} . \]   
For $l \leq k$ set
\[ \kappa_{k,l} = \iota_{k,l} \wedge \Psi :   bu \wedge (F_{4k}/F_{4k-1}) \wedge   \Sigma^{-2}  \mathbb{CP}^{2} \longrightarrow 
 bu \wedge (F_{4l}/F_{4l-1}) \wedge   \Sigma^{-2}  \mathbb{CP}^{2}  \]
then we obtain the following result.
\end{em}
\end{dummy}
\begin{theorem}{$_{}$}
\label{5.0b}

In the notation of \S\ref{5.0a}, in the $2$-local stable homotopy category, there exists $C' \in  {\rm Aut}_{left-bu-mod}^{0}( bu \wedge bo ) $
such that
\[  1 \wedge \psi^{3} : bu \wedge bu \longrightarrow  bu \wedge bu \]
satisfies
\[ ( C' \wedge 1) (1 \wedge \psi^{3}) (C' \wedge 1)^{-1} = \sum_{k \geq 0} \ 9^{k} \kappa_{k,k} + \sum_{k \geq 1} \ \kappa_{k,k-1} .   \]
\end{theorem}
\begin{dummy}{ $End_{left-bu-mod}(bu \wedge bo)$}
\label{5.1}
\begin{em}

In this section we shall apply Theorem \ref{0.1} to study the ring of left-$bu$-module homomorphisms of $bu \wedge bo$. As usual we
shall work in the $2$-local stable homotopy category. Let $\tilde{U}_{\infty}\mathbb{Z}_{2}$ denote the ring of upper triangular,
infinite matrices with coefficients in the $2$-adic integers. Therefore the group $U_{\infty}\mathbb{Z}_{2} $ is a subgroup of the
multiplicative group of units of $\tilde{U}_{\infty}\mathbb{Z}_{2}$. Choose a left-$bu$-module homotopy equivalence 
of the form
\[   \hat{L} :  \vee_{k \geq 0}  bu \wedge  (  F_{4k}/F_{4k-1})  
\stackrel{\simeq}{\longrightarrow }  bu  \wedge  bo,  \]
as in \S\ref{2.1}.
For any matrix $A \in
\tilde{U}_{\infty}\mathbb{Z}_{2}$ we may define a left-$bu$-module endomorphism of $bu \wedge bo$, denoted by $\lambda_{A}$,
by the formula
\[   \lambda_{A} =  \hat{L} \cdot (   \sum_{0 \leq l \leq k}  \  A_{l,k}  \iota_{k,l} )   \cdot    \hat{L}^{-1} . \]
Incidentally here and throughout this section we shall use the convention that a composition of maps starts with the right-hand map,
which is the {\em opposite} convention used in the definition of the isomorphism $\Psi$ of \S\ref{2.1} and \cite{SnTri}.
When $A \in U_{\infty}\mathbb{Z}_{2}$ we have the relation $\lambda_{A} = \Psi(A^{-1})$. For $A, B \in \tilde{U}_{\infty}\mathbb{Z}_{2}$
we have
\[ \begin{array}{ll}
\lambda_{A} \cdot \lambda_{B} & = (  \hat{L} \cdot (   \sum_{0 \leq l \leq k}  \  A_{l,k}  \iota_{k,l} )   \cdot   
\hat{L}^{-1} )
\cdot (  \hat{L} \cdot (   \sum_{0 \leq t \leq s}  \  B_{t,s}  \iota_{s,t} )   \cdot    \hat{L}^{-1} )  \\
\\
& =     \hat{L} \cdot (   \sum_{0 \leq l \leq t \leq s}  \  A_{l,t} 
 \  B_{t,s}     \iota_{s,l} )   \cdot    \hat{L}^{-1} \\
\\
& =    \hat{L} \cdot (   \sum_{0 \leq l  \leq s}  \  (AB)_{l,s} 
  \iota_{s,l} )   \cdot    \hat{L}^{-1} \\
\\
& = \lambda_{AB}.
\end{array}  \]

By Theorem \ref{0.1} there exists $H \in U_{\infty}\mathbb{Z}_{2}$ such that
\[ 1 \wedge \psi^{3} =  \lambda_{H B H^{-1}} \]
for 
\[  B = \left( 
Ê\begin{array}{cccccc} 

1 & 1 & 0 & 0 & 0 & \ldots \\
\\
0 & 9 & 1 & 0 & 0 & \ldots \\
\\
0 & 0 & 9^{2} & 1 & 0 Ê& \ldots \\
\\
0 & 0 & 0 & 9^{3} & 1 Ê & \ldots \\
\\
\vdots & Ê\vdots & Ê\vdots & Ê\vdots & Ê\vdots & Ê\vdots 
\end{array} Ê\right).  \]
Hence, for any integer $u \geq 1$, we have  $1 \wedge ( \psi^{3} - 9^{u-1}) =  \lambda_{HB_{u} H^{-1}}$
where $B_{u} = B - 9^{u-1} \in \tilde{U}_{\infty}\mathbb{Z}_{2}$ and $9^{u-1}$ denotes $9^{u-1}$ times the identity matrix.
Following \cite{Mg} write $\phi_{n} : bo \longrightarrow  bo$ for the composition
$\phi_{n} = (\psi^{3} - 1)(\psi^{3} - 9) \ldots  (\psi^{3} - 9^{n-1})$. Write $X_{n} = B_{1}B_{2} \ldots  B_{n} \in 
\tilde{U}_{\infty}\mathbb{Z}_{2}$.
\end{em}
\end{dummy}
\begin{theorem}
\label{5.2}

(i)  \  In the notation of \S\ref{5.1} $1 \wedge \phi_{n} = \lambda_{HX_{n}H^{-1}}$ for $n \geq 1$.

(ii)   \  The first $n$-columns of $X_{n}$ are trivial. 

(iii)  \   Let $C_{n} = Cone(  \hat{L} :  \vee_{0  \leq k  \leq n-1}  bu \wedge  (  F_{4k}/F_{4k-1})  
\stackrel{\simeq}{\longrightarrow }  bu  \wedge  bo,  $, which is a left-$bu$-module spectrum. Then in the 2-local stable
homotopy category there exists a commutative diagram of left-$bu$-module maps of the form
\newline
\begin{picture}(259,182)
\thinlines    
              \put(178,102){$C_{n}$}
              \put(293,166){$bu \wedge bo$}
              \put(66,166){$bu \wedge bo$}
              \put(89,160){\vector(3,-2){65}}
              \put(215,112){\vector(3,2){65}}
              \put(120,168){\vector(1,0){130}}
              \put(180,173){$1  \wedge \phi_{n}$}
              \put(100,133){$ \pi_{n}$}
              \put(270,133){$\hat{ \phi}_{n}$}
\end{picture}
\newline
where $\pi_{n}$ is the cofibre of the restriction of $\hat{L}$. Also $ \hat{\phi}_{n}$ is determined up 
to homotopy by this diagram.

(iv)  \  More precisely, for $n \geq 1$ we have

\[  (X_{n})_{s,s+j} = 0 \  {\rm if}  \  j<0  \  {\rm or}  \  j> n   \]
and the other entries are given by the formula
\[  (X_{n})_{s,s+t} =
\sum_{1 \leq k_{1} < k_{2} < \ldots < k_{t} \leq n}  \   A(k_{1}) A(k_{2}) \ldots  A(k_{t}) \]
where
\[ \begin{array}{l}
 A(k_{1}) = \prod_{j_{1} = n - k_{1} + 1}^{n}  (9^{s-1} - 9^{j_{1}-1}),  \\
\\
A(k_{2}) =  \prod_{j_{2} = n - k_{2} + 1}^{n- k_{1} -1}  (9^{s} - 9^{j_{2}-1}),  \\
\\
A(k_{3}) =  \prod_{j_{3} = n - k_{3} + 1}^{n- k_{2} -1}  (9^{s+1} - 9^{j_{3}-1}),  \\
\\
\hspace{30pt}  \vdots  \hspace{30pt}   \vdots \hspace{30pt}  \vdots   \hspace{30pt}   \vdots \\
\\
A(k_{t})  =    \prod_{j_{t} = 1}^{n- k_{t} -1}  (9^{s+ t + 1} - 9^{j_{t}-1})  .
\end{array}  \]

\end{theorem}
\vspace{2pt}

{\bf Proof}
\vspace{2pt}

Part (i) follows immediately from the discussion of \S\ref{5.1}. Part (ii) follows from part (iv), but it is
simpler to prove it directly. For part (ii) observe
that the $B_{i}$ commute, being polynomials in the matrix $B$ so that $X_{n} = X_{n-1} B_{n}$. Since $(B_{n})_{s,t}$ is 
zero except when $t=s, s+1$ so that $(X_{n})_{i,j}  = ( X_{n-1})_{i,j} (B_{n})_{j,j} + (X_{n-1})_{i,j-1} ( B_{n})_{j-1,j}$, which is
zero by induction if $j<n$. When $j=n$ by induction we have $(X_{n})_{i,j}  = (X_{n-1})_{i,n} (B_{n})_{n,n} $ which is trivial
because $(B_{n})_{n,n} = 9^{n-1} - 9^{n-1}$. In view of the decomposition of $bu \wedge bo$, part (iii) amounts to showing that
$HX_{n}H^{-1}$ corresponds to a left-$bu$-module endomorphismm of $\vee_{0 \leq k} \ bu \wedge (F_{4k}/F_{4k-1})$
which is trivial on each summand  $ bu \wedge (F_{4k}/F_{4k-1})$ with $k \leq n-1$. The $(i,j)$-th entry in this matrix is the 
multiple of $\iota_{j-1,i-1} : bu \wedge (F_{4j-4}/F_{4j-5}) \longrightarrow   bu \wedge (F_{4i-4}/F_{4i-5})$ given 
by the appropriate component of the map. The first $n$ columns are zero if and only if the map has no non-trivial components 
whose domain is $bu \wedge (F_{4j-4}/F_{4j-5})$ with $j \leq n$. Since $H$ is upper triangular and invertible, the first $n$ columns
of $X_{n}$ vanish if and only if the same is true for $HX_{n}H^{-1}$. Finally the formulae of part (iv) result from
the fact that $B_{j}$ has $9^{m-1} -  9^{j-1}$ in the $(m,m)$-th entry, $1$ in the $(m,m+1)$-th entry and zero elsewhere. $\Box$
\begin{remark}
\label{5.3}
\begin{em}

Theorem \ref{5.2} is closely related to the main result of \cite{Mg}. Following \cite{Mg} let 
$bo^{(n)} \longrightarrow  bo$ denote the map of $2$-local spectral which is universal for all maps $X \longrightarrow bo$ which are
trivial with respect to all higher $\mathbb{Z}/2$-cohomology operations of order less than $n$. Cf with
\cite{Mg} Theorem B. Milgram shows that $\phi_{2n}$ factorises through a map of the form 
$ \theta_{2n} : bo \longrightarrow  \Sigma^{8n}bo^{(2n - \alpha(n))}$ and that 
 $\phi_{2n+1}$ factorises through a map of the form 
$ \theta_{2n+1} : bo \longrightarrow  \Sigma^{8n+4} bsp^{(2n - \alpha(n))}$ and then uses the $\theta_{m}$'s to 
produce a left-$bo$-module splitting of $bo \wedge bo$.
Using the homotopy equivalence $bu \simeq bo \wedge  \Sigma^{-2} \mathbb{CP}^{2}$ mentioned in \cite{SnTri}
one may pass from the splitting of $bu \wedge bo$ to that of $bo \wedge bo$ (and back again). In the light of this observation, the
existence of the diagram of Theorem \ref{5.2} should be thought of as the upper triangular matrix version of the proof that the
$\theta_{n}$'s exist. The advantage of the matrix version is that Theorem \ref{5.2}(iv) gives us every entry in the matrix
$X_{n}$, not just the zeroes in the first $n$ columns.
\end{em}
\end{remark}

\end{em}
\end{dummy}

\begin{thebibliography}{Mg}
\bibitem{A}  J.F. Adams: {\em Stable Homotopy and Generalised Homology}; 
University of Chicago Press (1974).


\bibitem{BP} E.H. Brown and F.P. Peterson: On the stable decomposition of $\Omega^{2} S^{r+2}$;
Trans. Amer. Math. Soc. 243 (1978) 287-298.

\bibitem{ccw} F. Clarke, M. D. Crossley and S. Whitehouse: Bases for cooperations in $K$-theory;
$K$-Theory (3) 23 (2001) 237-250.
   


\bibitem{Mah} M. Mahowald: $bo$-Resolutions; Pac. J. Math. (2) 92 (1981) 365-383.

\bibitem{Mg} R.J. Milgram: The Steenrod algebra and its dual for 
connective K-theory; Reunion Sobre 
Teoria de Homotopia, Universidad de Northwestern 
Soc. Mat. Mex. (1974) 127-159.


\bibitem{Snspl} V.P. Snaith: A stable decomposition of $\Omega^{n} S^{n}X$; J. London Math. Soc.
2 (1974) 577-583.


\bibitem{SnTri}  V.P. Snaith: The upper triangular group and operations in algebraic K-theory;
Topology 41 (2002) 1259-1275.



\end{thebibliography}
\end{document}